\documentclass[11pt]{article}

\usepackage{epsfig}
\usepackage{latexsym}

\usepackage{epsfig}
\usepackage{latexsym}
\usepackage{amsmath,amssymb}
\usepackage{xcolor}
\input{epsf}

\newtheorem{prelem}{{\bf Theorem}}

\newenvironment{theoremABC}{\begin{prelem}{\hspace{-0.5
               em}{\bf.}}}{\end{prelem}}

\newtheorem{theorem}{Theorem}[section]

\newtheorem{definition}[theorem]{Definition}
\newtheorem{conjecture}[theorem]{Conjecture}

\newtheorem{lemma}[theorem]{Lemma}
\newtheorem{proposition}[theorem]{Proposition}

\newenvironment{proof}{{\bf Proof.}}{\hfill\rule{2mm}{2mm}}

\newenvironment{erratumblock}{\begingroup}{\endgroup}
\newtheorem{remarka}[theorem]{Remark}
\newenvironment{remark}{\begin{remarka}\rm}{\hfill\rule{2mm}{2mm}\end{remarka}}
\newtheorem{examplea}[theorem]{Example}
\newenvironment{example}{\begin{examplea}\rm}{\hfill\rule{2mm}{2mm}\end{examplea}}

\def\Ex {{\bf E}}
\def\Pr {{\rm Pr}}
\def\Var {{\bf Var}}

\def\DI {{I}}
\def\AI {{\tilde{I}}}

\title{\bf Decision trees and influences of variables over product probability spaces}
\author{Hamed Hatami
 \\
{\small\it Department of Computer Science }\\
{\small University of Toronto}\\ {\small\it e-mail:
hamed@cs.toronto.edu} }
\date{}
\begin{document}
\maketitle

\tableofcontents

\begin{abstract}
A celebrated theorem of Friedgut says that every function
$f:\{0,1\}^n \rightarrow \{0,1\}$ can be approximated by a function
 $g:\{0,1\}^n \rightarrow \{0,1\}$ with $\|f-g\|_2^2 \le \epsilon$ which depends only on
$e^{O(\DI_f/\epsilon)}$ variables where $\DI_f$ is the sum of the
influences of the variables of $f$. Dinur and Friedgut later showed
that this statement also holds if we replace the discrete domain
$\{0,1\}^n$ with the continuous domain $[0,1]^n$, under the extra
assumption that $f$ is increasing. They conjectured that the
condition of monotonicity is unnecessary and can be removed.

We show that certain constant-depth decision trees provide
counter-examples to Dinur-Friedgut conjecture. This suggests a
reformulation of the conjecture in which the function $g:[0,1]^n
\rightarrow \{0,1\}$ instead of depending on a small number of
variables has a decision tree of small depth. In fact we prove this
reformulation by showing that the depth of the decision tree of $g$
can be bounded by $e^{O(\DI_f/\epsilon^2)}$.

Furthermore we consider a second notion of the influence of a variable, and
study the functions that have bounded total influence in this sense. We
use a theorem of Bourgain to show that these functions have certain properties.
We also study the relation between the two different notions of influence.
\end{abstract}

\begin{erratumblock}
\section*{Erratum (August 2026)}
Theorem~3.4 is incorrect.  In its proof, the sentence
``On the other hand by (a generalization of) Theorem~C'' invokes a
generalization of the Bonami--Beckner inequality from the discrete cube
to an arbitrary product probability space.  The asserted generalization
is false, so the displayed inequality following that sentence is not
valid in this generality.  In fact, Theorem~3.4 itself is false.

For example, let $X=\{1,\ldots,r\}$ with the uniform probability
measure, let $n=2$, and let $f:X^2\to\{0,1\}$ be given by
$f(x_1,x_2)=1$ if and only if $x_1=x_2=1$.  Then
\[
\Var[f]=r^{-2}(1-r^{-2}),\qquad
\DI_f(1)=\DI_f(2)=r^{-1},
\]
and
\[
\AI_f(1)=\AI_f(2)=r^{-2}(1-r^{-1}).
\]
Consequently, the right-hand side of the inequality asserted in
Theorem~3.4 is
\[
\frac{20r^{-2}(1-r^{-1})}{\log r},
\]
which is smaller than $\Var[f]$ for all sufficiently large $r$.
Theorem~3.4 and its proof should therefore be disregarded.  No later
result in the paper uses Theorem~3.4.

 \end{erratumblock}

\noindent {{\sc AMS Subject Classification:} \quad 06E30 - 28A35}
\newline
{{\sc Keywords:} influence, Boolean function, decision tree,
junta.}

\section{Introduction \label{intro}}

The notion of the influence of a variable on a Boolean function over
a product probability space plays an important role in computer
science, combinatorics, statistical physics, economics and game
theory (see for example~\cite{KKL,BKKKL,Fr04,B1,F,LMN}). It is
usually the case that the functions whose variables satisfy certain
bounds on their influences have simple structures. Although results
of this type are studied extensively (see for
example~\cite{B1,F,H06,FKN}), still many basic questions are open (A
handful of them can be found in~\cite{Fr04}).

For a set $X$ and given a function $f: X^n \rightarrow \mathbb{R}$,
roughly speaking the influence of the $i$-th variable measures the
sensitivity of $f$ with respect to the changes in the $i$-th
coordinate. There are numerous ways to measure this sensitivity and
each gives rise to a different notion of the influence. It is
usually the case that certain bounds on the influences imply some
structure on $f$. For example there is a theorem of Friedgut
(Theorem~\ref{thm:Friedgut} below) in this flavor that concerns the
case $X=\{0,1\}$. This theorem has many applications and provides a
tool for proving hardness of approximation results using PCP
reductions, and also serves as a tool in machine
learning~\cite{DS2005,KR2003,DRS02,OS06}. Dinur and Friedgut later
conjectured that this result holds in the continuous setting of
$X=[0,1]$ as well. In this article we give a counter-example to
their conjecture. Fortunately this does not mean that there is no
continuous version of Friedgut's theorem. Indeed we characterize all
functions that satisfy the conditions of Dinur-Friedgut conjecture,
and prove that such functions enjoy very nice and simple structures.

A graph on $n$ vertices can be thought of as an element in
$\{0,1\}^{n \choose 2}$, where each coordinate corresponds to one of
the ${n \choose 2}$ possible edges. A graph property on $n$ vertices
is a subset of graphs on $n$ vertices which is invariant under the
permutation of vertices. Hence such a graph property can be
formulated as a function $f: \{0,1\}^{n \choose 2} \rightarrow
\{0,1\}$ which is invariant under certain permutations of
coordinates: the ones that come from the permutations of the
vertices of a graph. Let $\mu_p(f)$ denote the probability that
$f(G(n,p))=1$, where $G(n,p)$ is the random graph on $n$ vertices
where each edge is present with probability $p$. It follows from a
theorem of Margulis~\cite{Margulis} and Russo~\cite{Russo} that
$\frac{d \mu_p}{dp}$ is closely related to the influences of
variables for $f$. The derivative $\frac{d \mu_p}{dp}$ plays an
important role in theory of random graphs and statistical physics. A
bounded derivative shows that the property changes smoothly with
respect to the changes in $p$. Friedgut~\cite{FriBou} used the
connection to the influences to prove a beautiful theorem that
characterizes all graph properties that satisfy $\frac{d \mu_p}{dp}
< C$ for a constant $C$. He conjectured that being a graph property
is irrelevant here and his characterization holds for every function
$f:\{0,1\}^n \rightarrow \{0,1\}$. Bourgain~\cite{FriBou} proved a
weaker version of Friedgut's conjecture. We will observe that
Bourgain's proof works in the more general setting of functions
$f:X^n \rightarrow \mathbb{R}$. Inspired by this we give a
conjecture that is similar to Friedgut's conjecture but is for the
general setting of $X^n$. We show that a positive answer to this
conjecture would settle one of the conjectures in~\cite{Fr04}.

\section{Definitions and known facts}

\subsection{Notation}
Throughout this paper, following Hardy notation, $C$ always denotes
a universal constant, but not necessarily the same in different
statements. Let $[n]=\{1,\ldots, n\}$ for every natural number $n$.
For two functions $f,g: \mathbb{R}^+ \rightarrow \mathbb{R}^+$, we
write $f=O(g)$ if there exists a constant $C$ such that $f(x) \le C
g(x)$, for every $x \ge C$.  We also write $f=o(g)$ if $\lim_{x
\rightarrow \infty} f(x)/g(x)=0$. Sometimes we write $f=O_a(g)$ or
$f=o_a(g)$ to specify that the involved constants may depend on $a$.
When there is no subscript, the constants are universal. A function
is called \emph{Boolean} if its range is $\{0,1\}$.

Let $(X,{\cal F},\mu)$ be a probability space, and let $X^n$ be the
product space endowed with the product probability measure $\mu^n$.
Our goal is to study measurable functions $f: X^n \rightarrow
\mathbb{R}$. Sometimes we consider $f$ as a random variable and
write $\Ex[f]$ (the \emph{expected value} of $f$) for $\int f(x)
d\mu^n(x)$. The \emph{variance} of $f$ is defined as $\Var[f] =
\Ex[f^2] -(\Ex[f])^2$. For $J \subseteq [n]$, and $a \in X$, we
denote by $(J=a)$ the event $\{x: x_i=a \forall i \in J\}$.
Thus we may write $\Ex[f|J=a]$ to denote the conditional expectation
of $f$ on the event that all variables with indices in $J$ are set
to be equal to $a$.

Suppose that $X$ is an ordered set. Then a function $f:X^n
\rightarrow \mathbb{R}$ is called \emph{increasing} if $x_i \le y_i$
for all $i \in [n]$, implies $f(x) \le f(y)$.

\subsection{Fourier-Walsh expansion}
Let $(X,{\cal F},\mu)$ be a probability space, and let $X^n$ be the
product space endowed with the product probability measure $\mu^n$.
Consider a function $f:X^n \rightarrow \mathbb{R}$. Then it is
possible to find functions $F_S: X^n \rightarrow \mathbb{R}$ for $S
\subseteq [n]$ such that
\begin{itemize}
\item[{\bf (i)}] We have $f = \sum_{S \subseteq [n]} F_S$.

\item[{\bf (ii)}] For $x=(x_1,\ldots,x_n) \in X^n$, the value of $F_S(x)$ depends only on the variables $\{x_i: i
\in S\}$.

\item[{\bf (iii)}] We have $\int F_S d\mu(x_i)=0$ for every $i \in S$.

\item[{\bf (iv)}] We have $\int F_S F_T d\mu(x_i)=0$ for every $i \in T\setminus S$.
\end{itemize}
Then $f=\sum_{S \subseteq [n]} F_S$ is called the Fourier-Walsh
expansion of $f$. Note that (iv) is implied by (iii) and from (i-iv) it is easy to see that
\begin{equation}
\|f\|_2^2 = \int f^2 d\mu^n= \sum_{S \subseteq [n]} \sum_{T
\subseteq [n]} \int F_S F_T d\mu^n = \sum_{S \subseteq [n]} \int
F_S^2 d\mu^n = \sum_{S \subseteq [n]} \|F_S\|_2^2.
\end{equation}

For uniform measure on $X=\{-1,1\}$, the Fourier-Walsh expansion
becomes very simple. For $S \subseteq [n]$, let $w_S:X^n \rightarrow
\mathbb{R}$ be defined as $w_S: (x_1,\ldots,x_n) \mapsto \prod_{i
\in S} x_i$. Then it is easy to see that for $f:X^n \rightarrow
\mathbb{R}$,
$$f = \sum_{S \subseteq [n]} \widehat{f}(S) w_S,$$
where
$$\widehat{f}(S) = \int f(x) w_S(x)  d\mu^n(x).$$
So here in the Fourier-Walsh expansion $f=\sum_{S \subseteq [n]}
F_S$, we have $F_S=\widehat{f}(S) w_S$, where $\widehat{f}(S)$ is
just a real number.

\subsection{Influences\label{sec:influences}}
Let $(X,{\cal F},\mu)$ be a probability space, and let $X^n$ be the
product space endowed with the product probability measure $\mu^n$.
For any $x=(x_1,\ldots,x_n) \in X^n$, define
\begin{equation}
s_j(x)=\{y \in X^n:y_i = x_i \ \forall i \neq j\},
\end{equation}
or in other words $s_j(x)$ is the set of elements that can be
obtained from $x$ by changing only the $j$-th coordinate. For $f:X^n
\rightarrow \mathbb{R}$, we want to define the notion of the
influence of the $j$-th variable as a measurement for the dependence
of the value of $f$ to the value of the $j$-th coordinate.  Suppose
that $f$ is constant on $s_j(x)$ for some $x$. Then for this
particular $x$, the value of $x_j$ is irrelevant to the value of
$f$. Thus it is natural to define the influence of $x_j$ as the
probability that $f$ is not constant on $s_j(x)$. This gives rise to
our first notion of the influence. Next consider the function
$f:[0,1]^2 \rightarrow \mathbb{R}$ defined as
$f(x_1,x_2)=10x_1+x_2$. With the above definition both $x_1$ and
$x_2$ have influences $1$. But note that changing the value of $x_1$
changes the value of $f$ more drastically than changing the value of
$x_2$. The next definition of the influence captures this difference
and is defined as the following. For a particular $x$, $s_j(x)$ is
endowed with the probability measure $\mu$, and $\Var_{s_j(x)}[f]$
measures the variation of $f$ on $s_j(x)$. Thus we can define the
influence of $x_j$ as $\Ex [\Var_{s_j(x)}[f]]$. Let us summarize the
discussion in the following definition.

\begin{definition}
Let $f:X^n \rightarrow \mathbb{R}$ where $X^n$ is a product
probability space with the product probability measure $\mu^n$.
\begin{itemize}
\item The influence of the $i$-th variable on $f$ is defined as
\begin{equation}
\label{eq:digital} \DI_f(i)=\mbox{$\Pr(\{x: f$ is not constant on
$s_i(x)\})$}.
\end{equation}
\item The variational influence of the $i$-th variable on $f$ is defined
as
\begin{equation}
\label{eq:variational} \AI_f(i)=\Ex[\Var_{s_i(x)}[f]].
\end{equation}
\end{itemize}
The total influence and total variational influences on $f$ are
defined respectively as
$$
\mbox{$\DI_f=\sum_{i=1}^n \DI_f(i)$\qquad  and \qquad
$\AI_f=\sum_{i=1}^n \AI_f(i).$}
$$
\end{definition}

\begin{remark}
\label{rem:measurability} Note that (\ref{eq:digital}) and
(\ref{eq:variational}) imply some measurability conditions on $f$
which we assume to hold throughout this article. Since we are
dealing with probability spaces, for our purpose, without loss of
generality we might assume that $(X,\mu)$ is the Lebesgue measure on
$[0,1]$.  Throughout the article when a discrete set $X$ is
considered as a measure space without mentioning the measure, the
corresponding measure is assumed to be the uniform probability
measure.
\end{remark}

The influences are of particular interest when $f$ is the
characteristic function of a set, or in other words the range of $f$
is $\{0,1\}$.  Note that for a function $f:X^n \rightarrow
\{0,1\}$, we always have $4\AI_f(j) \le \DI_f(j)$. Indeed let
$$a_j(x)=\left\{\begin{array}{lcl}
               1 &\ &\mbox{$f$ is not constant on $s_j(x)$} \\
               0 & &\mbox{$f$ is a constant on $s_j(x)$}
         \end{array}\right..
$$
Let $x$ be a random variable that takes values in $X^n$ according to
$\mu^n$. Note that $a_j(x)=0$ implies $\Var_{s_j(x)}[f]=0$, and in
general since $f$ is Boolean, $\Var_{s_j(x)}[f] \le \frac{1}{4}$.
Thus
\begin{equation}
\label{eq:analog} \DI_f(j) = \Ex[a_j(x)] \ge 4\Ex [\Var_{s_j(x)}[f]]
=4\AI_f(j).
\end{equation}

In the discrete case of the uniform measure on $\{0,1\}^n$, the
equality holds in (\ref{eq:analog}), and the two definitions differ
only in a constant factor of $4$. But in the context of general
Boolean functions, $\DI_f(j)$ can be much larger than $\AI_f(j)$.
For example define $f:[0,1]^n \rightarrow \{0,1\}$ as
$f(x)=1$ if $x_j \ge \epsilon$, and $f(x)=0$ otherwise. Then
$\DI_f(j)=1$ while $\AI_f(j)=\epsilon(1-\epsilon)$.

The well-known result of~\cite{KKL} which is known as the KKL
inequality says that for every function $f:\{0,1\}^n \rightarrow
\{0,1\}$, there exists a variable whose influence is at least $C
\rho(1-\rho)\frac{\ln n}{n}$ where $\rho=\Ex[f]=\Pr_x[f(x)=1]$. Later
Friedgut~\cite{Fr04} noticed that the approach of \cite{KKL} leads
to an interesting result that later found many
applications~\cite{DS2005,KR2003,DRS02,OS06}:

\begin{theoremABC}\label{thm:Friedgut} For every $f:\{0,1\}^n \rightarrow \{0,1\}$, there exists
an approximation
 $g:\{0,1\}^n \rightarrow \{0,1\}$ with $\|f-g\|_2^2 \le \epsilon$ which depends only on
$e^{C\DI_f/\epsilon}$ variables.
\end{theoremABC}

In \cite{BKKKL}, Bourgain {\it et al.} generalized the framework of
KKL to the continuous domain, $f:[0,1]^n \rightarrow \{0,1\}$, and
proved that the KKL inequality holds in this setting as well. The
inequality in this case is usually called the BKKKL inequality.

Friedgut and Kalai~\cite{FK96} noticed that the proof
of~\cite{BKKKL} can be modified to imply the following statement
which is one of our main tools in this article.
\begin{theoremABC}\label{thm:BKKKL} For every $f: [0,1]^n \rightarrow \{0,1\}$ with
$\Ex[f] =\rho$, there always exists $i \in [n]$ such that
\begin{equation}
\label{eq:BKKKL}
 e^{-C \DI_f/\rho(1-\rho)} \le \DI_f(i).
 \end{equation}
\end{theoremABC}

The main tool in the proof of many results about the influences
including Theorem~\ref{thm:Friedgut} and Theorem~\ref{thm:BKKKL} is
the Bonami-Beckner inequality. This inequality was originally proved
by Bonami~\cite{Bonami} and then independently by
Beckner~\cite{Beckner}. It was first used to analyze discrete
problems in the proof of the KKL inequality~\cite{KKL}. We will use
the inequality in the following form.

\begin{theoremABC}[Bonami-Beckner Inequality]
\label{thm:bonami-beckner} Let $f:\{-1,1\}^n \rightarrow \mathbb{R}$
be a function with the Fourier-Walsh expansion $f:= \sum_{S \subseteq [n]}
\widehat{f}(S)w_S$. Then  for $p \ge 2$, and $\delta=\sqrt{p-1}$,
$$\|f\|_p\leq  \left\|\sum_{S \subseteq [n]} \delta^{|S|} \widehat{f}(S)w_S\right\|_2;$$
and for $1 \le p \le 2$,
$$\|f\|_p\ge \left\|\sum_{S \subseteq [n]} \delta^{|S|} \widehat{f}(S)w_S\right\|_2.$$
\end{theoremABC}

\subsection{Decision trees}

Consider a product probability space $X^n$ with the product
probability measure $\mu^n$. A \emph{decision tree} $T$ is a tree
where every internal node $v$ is labeled with an index $i_v \in [n]$
such that if $u$ is the ancestor of $v$, then $i_u \neq i_v$.
Children of $v$ are in a one-to-one correspondence with the elements
of $X$. Moreover a value $\mathrm{val}(v) \in \mathbb{R}$ is
assigned to every leaf $v$ of $T$. The \emph{depth} of the decision
tree is the maximum number of edges on a path from the root to a
leaf. We call the immediate parents of the leaves, the
\emph{leaf-parents}. A decision tree is called \emph{complete} if
all leaves are in the same distance from the root.

Every decision tree $T$ corresponds to a function $f_T:X^n
\rightarrow \mathbb{R}$: consider $x \in X^n$; we start from the
root and explore a path to one of the leaves in the following way:
Every time that we meet an internal node $v$, we look at the value
$x_{i_v}$ and move to the corresponding child of $v$. We continue
until we reach a leaf $u$. Then $f(x)=\mathrm{val}(u)$. The
traversed path from the root to the leaf is called the
\emph{computing path} of $x$.

There are interesting relations between the structure of a decision
tree $T$, and the influences of the variables on $f_T$ (see for
example~\cite{Ryan05}). The following simple lemma shows that the
total influence on $f_T$ is bounded by the depth of $T$.

\begin{lemma}
\label{lem:treeinfluence} Consider a decision tree $T$ of depth
$k$ and the function $f_T$. Then $\DI_{f_T} \le k.$
\end{lemma}
\begin{proof}
Choose $x \in X^n$ at random according to $\mu^n$, and $i \in [n]$
uniformly at random. If $i$ does not appear on the computing path of
$x$, then $f_T$ is constant on $s_i(x)$. The probability that $i$
appears on the computing path is at most $k/n$. Now the assertion of
the lemma follows from the linearity of the expectation and the fact
that there are $n$ different variables.
\end{proof}

\section{Main Results}
\subsection{Counter-example to Dinur-Friedgut Conjecture \label{sec:counterexample}}

The proof of BKKKL inequality and Theorem~\ref{thm:BKKKL} rely on a
simple reduction that reduces the problem to the case that $f$ is
increasing. Similar but more complicated arguments of this sort were
previously used by Talagrand in the context of isoperimetric
inequalities for Banach-Space valued random variable~\cite{Tal89}.

A key but simple observation in \cite{BKKKL} is that in the case of
increasing functions, one can without increasing the influences by a
factor of more than 2, reduce the problem further to the case where
the domain of the function is $\{0,1\}^m$ endowed with the uniform
probability measure, where $m=O(n\log n)$.

Although Theorem~\ref{thm:BKKKL} provides a continuous version of
the KKL inequality, a continuous version of
Theorem~\ref{thm:Friedgut} remained unknown. In fact Dinur and
Friedgut (see \cite{Fr04}) considered the case where $f$ is
increasing, and observed that in this case the reduction of
\cite{BKKKL} to $\{0,1\}^m$ shows that $f$ essentially depends on
$e^{CI_f/\epsilon}$ number of variables, or formally for every
$\epsilon>0$ there exists a function $g:[0,1]^n \rightarrow \{0,1\}$
which depends on $e^{CI_f/\epsilon}$ number of variables and
$\|f-g\|_2^2 \le \epsilon$. This provides a continuous version of
Theorem~\ref{thm:Friedgut}, but under the extra assumption of
monotonicity. They conjectured (Conjecture~2.13 in \cite{Fr04}) that
this extra condition is unnecessary. The reason that the reduction
to increasing functions cannot be used here is that with this
reduction one may get a function which depends on a small number of
variables while the original function cannot be approximated by any
such function.

We disprove Dinur-Friedgut conjecture in the following lemma by
introducing a function for which the total influence is small, but
it cannot be approximated by any function which depends on a small
number of variables:

\begin{lemma}
\label{lem:counterexample} Let $k \in \mathbb{N}$ and $r \ge 2$ be
an even integer, and let $\mu$ be the uniform probability measure on
$X=\{0,\ldots,r-1\}$. For $n \ge r^k$, there exists a function
$f:X^n \rightarrow \{0,1\}$ such that $\DI_f \le k$ and $\|f-g\|_2^2
\ge \frac{1}{4}$ for every $g:X^n \rightarrow \{0,1\}$ that depends
on at most $r^{k-1}/2$ variables.
\end{lemma}
\begin{proof}
Consider a complete decision tree $T$ of depth $k$ where $i_v \neq
i_u$ for every two internal nodes $u$ and $v$. Consider a leaf $w$
which is the immediate child of a node $u$, and is corresponded to
$a \in \{0,\ldots,r-1\}$. We let $\mathrm{val}(w)=\lfloor
\frac{2a}{r} \rfloor$.

By Lemma~\ref{lem:treeinfluence} we have $\DI_{f_T} \le k$. On the
other hand consider a function $g:X^n \rightarrow \{0,1\}$ that
depends on at most $r^{k-1}/2$ variables. Let $L$ be the set of the
leaf-parents of $T$. Note that $|L| = r^{k-1}$ which shows that for
at least half of the nodes $v \in L$, $g$ is a constant on the
children of $v$. This completes the proof because $f_T \neq g$ on
half of the children of every such $v$.
\end{proof}

In the first look Lemma~\ref{lem:counterexample}  may sound
disappointing because it shows that the continuous version of
Theorem~\ref{thm:Friedgut} is not true. But further inspection shows
that indeed a continuous version does exist, but it is slightly
different from the discrete case. Our counter-example to
Dinur-Friedgut conjecture has a decision tree of constant depth. We
will show that for a fixed $\epsilon$, having a constant-depth
decision tree is the right reformulation of the conjecture. First we
``discretize'' the problem:

\begin{remark}\label{rem:discrete}
Let $f:[0,1]^n \rightarrow \{0,1\}$ be a function that satisfies the
measurability conditions of Remark~\ref{rem:measurability}. Fix any
$0<\delta<\frac{1}{10}$. Consider $m \in \mathbb{N}$ and subdivide
$[0,1]^n$ into $m^n$ equal size disjoint cells by subdividing each
one of the base intervals into $m$ intervals. For every $x \in
[0,1]^n$ denote by $C_x$ the cell that contains $x$. For every $1
\le i \le n$, and cell $C$, let $s_i(C)$ denote the set of $m$ cells
that are obtained by changing the cell $C$ in the $i$th coordinate.
Consider a measurable set $A \subseteq [0,1]^n$. We say that a cell
$C$ is $(1-\delta,A)$-\emph{determined}, if either $\Pr(A | C) \ge
1-\delta$ or $\Pr(A^c | C) \ge 1-\delta$. Since $A$ is measurable,
for sufficiently large $m$, more than $(1-\delta)$ fraction of the
cells will be $(1-\delta,A)$-determined. For every $1 \le i \le n$,
let
$$A_i:=\{x: \mbox{$f$ is a constant on $s_i(x)$}\}.$$
Since $A_i$ are measurable, one can take $m$ to be sufficiently
large so that $(1-\delta)$ fraction of the cells are simultaneously
$(1-\delta,f)$-determined and $(1-\delta,A_i)$-determined, for every
$1 \le i \le n$.

Now we construct a ``discrete'' approximation $h$ of $f$, in the
following way:
\begin{enumerate}
\item[(i)] For every cell $C$, $h(C)=\{1\}$, if $\Ex(f|C) \ge 1-\delta$. Furthermore for every such $C$ and every $1 \le i \le n$,
 if $\Pr(A_i | C) \ge
1-\delta$, then $h(C')=\{1\}$ for every $C' \in s_i(C)$.

\item[(ii)] Similarly for every cell $C$, $h(C)=\{0\}$, if $\Ex(1-f|C) \ge 1-\delta$.
Furthermore for every such $C$ and every $1 \le i \le n$, if
$\Pr(A_i | C) \ge 1-\delta$,  then $h(C')=\{0\}$ for every $C' \in
s_i(C)$.

\item[(iii)] The value of $h$ on the rest of the cells is defined
arbitrarily.
\end{enumerate}
First we need to verify that $h$ is well-defined. Note that for
every $C'$ in (i) we have $\Ex(f|C') \ge 1-2\delta>1/2$, while for
every $C'$ in (ii), $\Ex(1-f|C') \ge 1-2\delta>1/2$. This shows that
$h$ is well-defined.

Note that $h$ is defined so that if $C_x$ is
$(1-\delta,f)$-determined and $\Pr(A_i | C_x) \ge 1-\delta$, then
$h$ is constant on $s_i(C_x)$. Thus
\begin{eqnarray}
\nonumber \DI_h(i) &\le& \Pr(\{x : \mbox{$C_x$  is not
$(1-\delta,f)$-determined}\}) + \Pr(\{x : \mbox{$C_x$  is not
$(1-\delta,A_i)$-determined}\}) \\& & + \Pr(\{x : \Pr(A_i^c | C_x)
\ge 1-\delta\}) \le  \delta+\delta+ \frac{\Pr(A_i^c)}{1-\delta} =
2\delta+ \frac{\DI_f(i)}{1-\delta}. \label{eq:approxInf}
\end{eqnarray}
Moreover
\begin{eqnarray}
\nonumber  \Pr[f \neq h] &\le& \Pr(\{x : \mbox{$C_x$  is not
$(1-\delta,f)$-determined}\}) + \delta \Pr(\{x : \mbox{$C_x$  is
$(1-\delta,f)$-determined}\}) \\ &\le& 2\delta.
\label{eq:approxFunc}
\end{eqnarray}
Inequalities (\ref{eq:approxInf}) and (\ref{eq:approxFunc}) show
that by taking $\delta$ to be arbitrarily small, we can obtain a
``discrete'' approximation $h$ of $f$, which is arbitrarily close to
$f$ in $\ell_2$ norm, and its total influence is not much larger
than $\DI_f$.
\end{remark}

\begin{theorem}
\label{thm:decision} Let $f:[0,1]^n \rightarrow \{0,1\}$ satisfy
$\DI_f \le B$. Then for every $\epsilon>0$ there exists a function
$g:[0,1]^n \rightarrow \{0,1\}$ such that $\|f-g\|_2^2 \le \epsilon$
and $g$ has a decision tree of depth at most $e^{C B/\epsilon^2}$.
\end{theorem}
\begin{proof}
By Remark~\ref{rem:discrete}  we can assume that the underlying
probability space is the uniform probability measure on
$\{0,\ldots,r-1\}^n$ for a sufficiently large $r$ that depends on
$f$ and $\epsilon$.

By a \emph{raw decision tree} we mean the mathematical object that
is obtained by removing the values $\mathrm{val}(v)$ from the leaves
of a decision tree. Consider a function $f:\{0,\ldots,r-1\}^n
\rightarrow \{0,1\}$ with $\DI_f \le B$. For every node $v$ in a raw
decision tree, $f$ induces a function $f_v$ of the variables that
are not assigned a value through the path from the root to $v$. This
function is obtained by fixing the values of the variables on the
path from the root to $v$. Denote by $m(v)$ the variable with
maximum influence in $f_v$, and note that if $v$ is in distance $k$
from the root then $\mu(v)=r^{-k}$ is the fraction of $x\in
\{0,\ldots,r-1\}^n$ that $v$ belongs to their computing paths (i.e.
$x$ agrees with the path from the root to $v$).

Let $C' \ge 1$\footnote{We can assume that $C' \ge 1$ as
(\ref{eq:BKKKL}) remains valid if one increases the constant in the
exponent.} be the constant from the exponent in (\ref{eq:BKKKL}).
Start with a single node as a raw decision tree. Consider the
following procedure: For every leaf $v$, we let $i_v:=m(v)$ and we
branch on this variable, adding $r$ children to $v$. Each
application of this procedure increases the depth of the tree by
one. We repeatedly apply the procedure until we obtain a raw
decision tree $T$ in which the fraction of the leaf-parents $v$ that
satisfy $I_{f_v}(i_v) \ge e^{-18C' B/\epsilon^2}$ is at most
$\epsilon/3$. Let $S$ denote the set of the nodes $v$  that satisfy
$I_{f_v}(i_v) \ge e^{-18C' B/\epsilon^2}$, and let $L_1 \subseteq S$
be the set of the leaf-parents in $S$. We claim that the depth of
$T$ is at most $e^{(18C'+4)B/\epsilon^2}$.

Suppose to the contrary that this is not true. Consider a node $u$
and note that at node $u$, the $i_u$'th variable contributes $\mu(u)
\times \DI_{f_u}(i_u)$ to $\DI_f(i_u)$. Now consider another node $v
\in T$ such that $i_v=i_u$. Since none of $u$ and $v$ is the
ancestor of the other, the contribution of the $i_u$'th variable to
$I_f(i_u)$ at node $u$ is disjoint from its contribution at node
$v$. Denoting by $h$ the depth of $T$ and by $D_k$ the set of the
vertices in distance $k$ from its root, the above discussion shows
that
\begin{eqnarray*}
\DI_f &\ge& \sum_{v \in T} \mu(v) I_{f_v}(i_v) = \sum_{k=0}^{h-1}
\sum_{v \in D_k}  \mu(v) I_{f_v}(i_v) \ge \\ &\ge& \sum_{k=0}^{h-1}
\sum_{v \in D_k \cap S}  \mu(v) I_{f_v}(i_v) \ge \sum_{k=0}^{h-2}
(\epsilon/3) e^{-18C' B/\epsilon^2}.
\end{eqnarray*}
Thus $h \le (3B/\epsilon) e^{18C' B/\epsilon^2}+1 \le
e^{(18C'+4)B/\epsilon^2}$.

It is easy to see that for a similar reason, $\sum_{v} \mu(v)
\DI_{f_v} \le \DI_f$ where now the sum is only over all leaf-parents
of $T$. This shows that the fraction of the leaf-parents $v$ that
satisfy $\DI_{f_v} \ge 3B/\epsilon$ is at most $\epsilon /3$. Let
$L_2$ denote the set of the leaf-parents that satisfy $\DI_{f_v} \ge
3B/\epsilon$.

Consider a leaf-parent $v \not\in L_1 \cup L_2$, i.e.
$I_{f_v}(m(f_v)) \le e^{-18 C' B/\epsilon^2}$ and $I_{f_v} \le
3B/\epsilon$. By Theorem~\ref{thm:BKKKL},
$$\left(\Ex[f_v]\right)\left(1-\Ex[f_v]\right) \le \epsilon/6,$$
which implies that $\min(\Ex[f_v],1-\Ex[f_v]) \le \epsilon/3$. Now
to every leaf $w$ with parent $u$ we assign
$$
\mathrm{val}(w)= \left\{ \begin{array}{lcc} 0 & &  \Ex[f_u] \le 1/2 \\
1 & & \Ex[f_u] > 1/2 \end{array}\right..
$$
Let $g$ be the function that is computed by this decision tree.
Denote by $L$ the set of all leaf-parents and correspond to every
$x$ the unique leaf-parent $l(x) \in L$ on its computing path.
\begin{eqnarray*}
 \|f-g\|_2^2 &=& \Pr[f(x) \neq g(x)] \le \Pr[l(x) \in L_1]+\Pr[l(x) \in L_2]
 +\\
&&\sum_{v \in L \setminus (L_1 \cup L_2)}
\mu(v)\min(\Ex[f_v],1-\Ex[f_v]) \le
\frac{\epsilon}{3}+\frac{\epsilon}{3}+\frac{\epsilon}{3} \le
\epsilon.
\end{eqnarray*}
\end{proof}

\subsection{Variational Influences for Boolean functions}

Our next goal is to study the relation between the influences and
the variational influences. The next theorem shows that when all
influences are small then the total variational influence is large.
One can prove this theorem easily using Theorem~1.5
in~\cite{Talagrand94}, but for the sake of the completeness we give
a direct proof here.

\medskip
{\color{red}\noindent\textbf{\textcolor{red}{Erratum.}} The following theorem (Theorem~3.4) is false; see the erratum at the beginning of the paper.}
\medskip

\begin{theorem}
For $f:X^n \rightarrow \mathbb{R}$, we have
\begin{equation}
\label{eq:thmStKKL} \Var[f] \le 10 \sum_{i=1}^n \frac{\AI_f(i)
}{\log(1/\DI_f(i))}.
\end{equation}
\end{theorem}
\begin{proof}
Let $\iota_i(x)=f(x)-\Ex_{s_i(x)}[f]$, and note that $\DI_f(i) \ge
\Ex [1_{[\iota_i \neq 0]}]$,  and $\AI_f(i)=\|\iota_i\|_2^2$. Let
$f=\sum_{S \subseteq [n]} F_S$ be the Fourier-Walsh expansion of
$f$. Note that
$$\Ex_{s_i(x)}[F_S]= \int F_S d\mu(x_i) = \left\{\begin{array}{lc} F_S & i \not\in S \\ 0 & i \in S\end{array}\right.$$
Thus the Fourier-Walsh expansion of $\iota_i$ is the following.
$$\iota_i = \sum_{S
\subseteq [n]} F_S - \Ex_{s_i(x)}\left[\sum_{S \subseteq [n]} F_S\right]=
\sum_{S \subseteq [n]} F_S - \sum_{S: i \not\in S} F_S=\sum_{S:i \in
S} F_S.$$
For $1<p \le 2$, by H\"{o}lder's inequality we have
$$\|\iota_i\|_p=\left(\int \iota_i^p \right)^{1/p} =
\left(\int \iota_i^p 1_{[\iota_i \neq 0]} \right)^{1/p} \le
\|\iota_i\|_2 \left(\int 1_{[\iota_i \neq 0]}\right)^{(1-p/2)/p} \le
\|\iota_i\|_2 \DI_f(i)^{\frac{1}{p}-\frac{1}{2}}.$$
On the other hand by (a generalization of)
Theorem~\ref{thm:bonami-beckner}, we have
$$\left(\sum_{S:i \in S, |S| = k} \|F_S\|_2^2  \right)^{1/2} = \left(\sum_{|S| = k} \widehat{(\iota_i)}_S^2  \right)^{1/2} \le
(p-1)^{-k/2} \|\iota_i\|_p.$$
Thus for $p=3/2$ we get
$$\sum_{S: i \in S, |S| \le k} \frac{\|F_S\|_2^2}{|S|} =
\sum_{t=1}^k \frac{1}{t}\sum_{S: i \in S, |S| = t} \|F_S\|_2^2 \le
\sum_{t=1}^k \frac{2^t}{t} \AI_f(i)  \DI_f(i)^{1/3} \le
\frac{2^{k+2}}{k+1} \AI_f(i)  \DI_f(i)^{1/3}.$$
On the other hand since  $\sum_{S: i \in S} \|F_S\|_2^2 = \AI_f(i)$,
we get
$$\sum_{S: i \in S, |S| > k} \frac{\|F_S\|_2^2}{|S|} \le \frac{\AI_f(i)}{k+1},$$
and hence
$$\sum_{S: i \in S} \frac{\|F_S\|_2^2}{|S|} \le
\frac{1}{k+1}\left(2^{k+2}\AI_f(i)  \DI_f(i)^{1/3} +
\AI_f(i)\right).$$
Letting $k=\frac{\log 1/\DI_f(i)}{3}$ we get
\begin{equation}
\label{eq:boundEnt} \sum_{S: i \in S} \frac{\|F_S\|_2^2}{|S|} \le
\frac{10 \AI_f(i)}{\log 1/\DI_f(i)}.
\end{equation}
Summing (\ref{eq:boundEnt}) over all $i \in [n]$, we
get~(\ref{eq:thmStKKL}).
\end{proof}

 Let us next investigate the properties of the functions that satisfy $\AI_f \le B$ for some
constant $B$. As it is already noticed by Dinur and Friedgut
(see~\cite{Fr04}) a hypercontractivity argument similar to~\cite{F}
implies that every function $f:\{0,1,\ldots,r-1\}^n \rightarrow
\{0,1\}$ with $\DI_f \le B$ essentially depends on $r^{CB
/\epsilon}$ number of variables. In fact such an argument implies
more:

\begin{lemma}
\label{lem:domianR} Let $f:\{0,1,\ldots,r-1\}^n \rightarrow \{0,1\}$
satisfy $\AI_f \le B$. Then for every $\epsilon>0$, there exists a
function $g:\{0,1,\ldots,r-1\}^n \rightarrow \{0,1\}$ which depends
on $r^{CB /\epsilon}$ number of variables, and $\|f-g\|_2^2 \le
\epsilon$.
\end{lemma}

Lemma~\ref{lem:counterexample} shows that Lemma~\ref{lem:domianR} is
sharp up to the constant in the exponent. In the more general
setting of functions $f:[0,1]^n \rightarrow \{0,1\}$, bounding the
total variational influence implies less structure on the function.
The following example shows that Theorem~\ref{thm:decision} does not
need to hold when one replaces the bound $\DI_f \le B$ with $\AI_f
\le B$.

\begin{example}
\label{example:example1} Let $f:[0,1]^n \rightarrow \{0,1\}$ be
defined as $f(x)=0$ if and only if $x \le
(1-\frac{1}{n},\ldots,1-\frac{1}{n})$. It is easy to see that
$\Ex[f] = 1-(1-1/n)^n \approx 1-1/e$ and $\AI_f \le 2$. Note that $f$ cannot be
approximated by a function with a decision tree of depth $o(n)$.
\end{example}

The following theorem which is implied by Bourgain's proof for
Proposition~1 in the appendix of~\cite{FriBou} shows that  bounded
variational influence implies some weak structure on the function.
\medskip
{\color{red}\noindent\textbf{Correction.} In the original version, the hypothesis $\Ex[f]=1/2$ was inadvertently omitted from the following theorem.}
\medskip

\begin{theorem}
\label{thm:Bourgain} Let $f:[0,1]^n \rightarrow \{0,1\}$ satisfy
$\Ex[f]=1/2$ and $\AI_f \le B$. Then there exists a set $J \subseteq [n]$ of size at
most $10B$ and an assignment of values to the variables in $J$ such
that conditioning on this partial assignment, the expected value of
$f$ changes by at least $3^{-500 B^2}$.
\end{theorem}

To see that why Bourgain's proof implies Theorem~\ref{thm:Bourgain} note that according to Equation~(2.18) in
the appendix of \cite{FriBou} if $\Ex[f]=1/2$ and $\AI_f \le B$, then
\begin{equation}
\label{eq:Bourgain} \Ex\left[\max_{J: |J| \le 10 B}|\Ex_J[f](x)
-\Ex[f]|\right] \ge 3^{-500 B^2}.
\end{equation}
This immediately implies Theorem~\ref{thm:Bourgain}. Inspired by Conjecture~1.5 in
\cite{FriBou}, we conjecture that Theorem~\ref{thm:Bourgain} can be
improved to the following:

\begin{conjecture}
\label{conj:StrongFri} There exists a function $k:\mathbb{R}^2
\rightarrow \mathbb{R}$ such that the following holds. Let
$B,\epsilon>0$ be two constants. For every $f:[0,1]^n \rightarrow
\{0,1\}$ that satisfies $\AI_f \le B$, there exists a set $J
\subseteq [n]$ of size $k(B,\epsilon)$ and an assignment of values
to the variables in $J$ such that conditioning on this partial
assignment, the expected value of $f$ is either less than $\epsilon$
or greater than $1-\epsilon$.
\end{conjecture}

\subsection{Subcubes of increasing sets}
The following theorem proved in~\cite{KKL} is a straightforward
corollary of the KKL inequality.

\begin{theoremABC}
\label{thm:corKKL} There exists a function $k:\mathbb{R}^2
\rightarrow \mathbb{R}^+$ such that the following holds. For every
increasing function $f:\{0,1\}^n \rightarrow \{0,1\}$ with
$\Var[f]=\rho$ and every $\epsilon>0$, there exists a set of
coordinates $J \subseteq [n]$ such that $|J| \le k(\rho,\epsilon)n/\log
n$ and $\Ex[f|J=1] \ge 1-\epsilon$.
\end{theoremABC}

Although it is wrongfully claimed in \cite{BKKKL} that one can use
Theorem~\ref{thm:BKKKL} to prove that the same statement holds for
functions $f:[0,1]^n \rightarrow \{0,1\}$, as it is noticed by
Friedgut in~\cite{Fr04}, the function in
Example~\ref{example:example1} shows that this is not true.

Friedgut~\cite{Fr04} suggested the following conjecture as the
continuous version of Theorem~\ref{thm:corKKL}.
\begin{conjecture}
\label{conj:Friedgut} There exists a function $k:\mathbb{R}^2 \rightarrow \mathbb{R}$ such that for every $\epsilon>0$, $\lim_{n \rightarrow \infty}\frac{k(\epsilon,n)}{n}=0$, and
the following holds. Let $f:[0,1]^n \rightarrow \{0,1\}$ be an
increasing function. There exists a set of
coordinates $J \subseteq [n]$ such that $|J| \le k(\epsilon,n)$ and either
$\Ex[f(x) | J=1] \ge 1-\epsilon$ or $\Ex[f(x) | J=0] \le \epsilon$.
\end{conjecture}

The next proposition shows that Conjecture~\ref{conj:StrongFri}
implies Conjecture~\ref{conj:Friedgut}.

\begin{proposition}
\label{pro:twoConjs}  If Conjecture~\ref{conj:StrongFri} is true,
then there exists a function $k:\mathbb{R}^2 \rightarrow \mathbb{R}$
such that for every $\epsilon>0$, $\lim_{n \rightarrow \infty}\frac{k(\epsilon,n)}{n}=0$, and
the following holds. Let $f:[0,1]^n \rightarrow \{0,1\}$
be an increasing function. Then there exists a set of coordinates $J
\subseteq [n]$ such that $|J| \le k(\epsilon,n)$ and either $\Ex[f(x)
| J=1] \ge 1-\epsilon$ or $\Ex[f(x) | J=0] \le \epsilon$.
\end{proposition}
\begin{proof}
For a set $J \subseteq [n]$, let $f|_{J=1}:[0,1]^n \rightarrow
\{0,1\}$ be defined as $f|_{J=1}: x \mapsto f(y)$ where
$$y_i=\left\{\begin{array}{lcl}1 &\qquad & i \in J \\ x_i & & i
\not\in J \end{array}\right.$$

Fix some $B>0$, and let $J_0:=\emptyset$. For $i \ge 0$, inductively
obtain $J_{i+1}$ from  $J_i$ as in the following. We will use the
notation $f_i=f|_{J_i=1}$.

If $\AI_{f_i} \ge B$, then there exists a variable $j \in [n]$ such
that $\AI_{f_i}(j) \ge B/n$. Let $J_{i+1}=J_i \cup \{j\}$. Note that
for every $y \in A_i:=\{x: \Ex_{s_j(x)}[f_i] \neq 0\}$, we have
$f_{i+1}(y)=1$. Thus $\Ex[f_{i+1}] \ge \Ex[1_{[x \in A_i]}]$ and
furthermore we have $\Ex[f_i] = \Ex[\Ex_{s_j(x)}[f_i]]$ . Hence
\begin{eqnarray*}
\Ex [f_{i+1}]-\Ex[f_i]  &\ge& \Ex[1_{[x \in A_i]}] -
\Ex[\Ex_{s_j(x)}[f_i]] = \Ex[1_{[x \in A_i]} - \Ex_{s_j(x)}[
f_i]]\\
&\ge&
\Ex[(\Ex_{s_j(x)}[f_i])(1-\Ex_{s_j(x)}[f_i])]=\Ex[\Var_{s_j(x)}[f_{i}]]
= \AI_{f_i}(j) \ge B/n.
\end{eqnarray*}

Continuing in this manner  we find a set $J_k$ for some $k \le n/B$
such that $|J_k| \le n/B$ and either $\Ex[f|J_k=1] \ge 1-\epsilon$ or
$\AI_{f_k} \le B$. But if $\AI_{f_k} \le B$, then assuming
Conjecture~\ref{conj:StrongFri} we can find a set $J'$ of size
$\kappa(B,\epsilon)$ such that $\Ex[f_k|J'=0] \le \epsilon$ or
$\Ex[f_k|J'=1] \ge 1-\epsilon$. If $\Ex[f_k|J'=1] \ge 1-\epsilon$, then for
$J=J_k\cup J'$ we have $\Ex[f|J=1]\ge 1-\epsilon$.  If instead
$\Ex[f_k|J'=0]\le\epsilon$, then, by monotonicity of $f$,
\[
\Ex[f|J_k=0,J'=0]\le \Ex[f|J_k=1,J'=0]
=\Ex[f_k|J'=0]\le\epsilon.
\]
Thus in either case there is a set $J$ of size at most
$n/B+\kappa(B,\epsilon)$ satisfying the desired conclusion.  For
fixed $\epsilon$, define
\[
k(\epsilon,n)=\inf_{B>0}\left(\frac{n}{B}+\kappa(B,\epsilon)\right).
\]
For every fixed $B>0$,
\[
\limsup_{n\to\infty}\frac{k(\epsilon,n)}{n}\le \frac1B,
\]
and since $B$ is arbitrary, $k(\epsilon,n)=o(n)$.
\end{proof}

\section*{Acknowledgements} The author wishes to thank Ehud Friedgut
for his valuable discussions.

\bibliographystyle{alpha}
\bibliography{influences}

\end{document}